\documentclass[number,citesort,seceqn,dvips]{arxbj}
\usepackage{upgreek}
\usepackage{mathrsfs}

\aid{0}
\volume{17}
\issue{4}
\pubyear{2011}
\firstpage{1136}
\lastpage{1158}
\doi{10.3150/10-BEJ308}

\makeatletter

\newtheorem{tthm}{Theorem}[section]
\newtheorem{cor}[tthm]{Corollary}
\newtheorem{lem}[tthm]{Lemma}
\newtheorem{prp}[tthm]{Proposition}

\newremark{remm}{Remark}[section]
\newcommand{\scr}[1]{\mathscr #1}

\newcommand{\R}{\mathbb{R}}
\newcommand{\ff}{\frac}
\newcommand{\sss}{\sqrt}
\newcommand{\dd}{\delta}

\newcommand{\N}{\mathbb{N}}

\newcommand{\DD}{\Delta}
\newcommand{\vv}{\varepsilon}
\newcommand{\rr}{\rho}
\newcommand{\GG}{\Gamma}
\newcommand{\ggam}{\gamma}
\newcommand{\nn}{\nabla}
\newcommand{\dsd}{\mathrm{d}}
\newcommand{\aal}{\alpha}
\newcommand{\si}{\sigma}
\newcommand{\F}{\scr F}

\newcommand{\e}{\mathrm{e}}
\newcommand{\ttil}{\tilde}

\newcommand{\ppbl}{\mathbb{P}}

\newcommand{\Z}{\mathbb{Z}}
\newcommand{\llam}{\lambda}

\newcommand{\E}{\mathbb{E}}

\newcommand{\LL}{\Lambda}
\newcommand{\Rank}{\operatorname{Rank}}
\newcommand{\B}{\mathbf{B}}

\makeatother

\begin{document}
\begin{frontmatter}

\title{Coupling for Ornstein--Uhlenbeck processes with jumps}

\runtitle{Ornstein--Uhlenbeck processes with jumps}

\begin{aug}
\author{\fnms{Feng-Yu} \snm{Wang}\ead[label=e1]{wangfy@bnu.edu.cn}\ead[label=e2]{F.Y.Wang@swansea.ac.uk}\corref{}}

\runauthor{F.-Y. Wang}
\address{School of Mathematical Sciences,
Beijing Normal
University, Beijing 100875, China and Department of Mathematics,
Swansea University, Singleton Park, SA2 8PP, UK\\ \printead{e1,e2}}
\end{aug}

\received{\smonth{3} \syear{2010}}
\revised{\smonth{6} \syear{2010}}

\begin{abstract}
Consider the linear stochastic differential equation (SDE) on $\mathbb{R}^n$:
\[
\dsd{X}_t= A X_t\,\mathrm{d}t+ B\,\mathrm{d}L_t,
\]
where $A$ is a real $n\times n$   matrix, $B$ is a real $n\times d$
real matrix and $L_t$ is a  L\'{e}vy  process  with L\'{e}vy measure
$\nu$ on $\mathbb{R}^d$.  Assume that $\nu(\dsd{z})\ge \rho_0(z)\,\mathrm{d}z$ for
some $\rho_0\ge 0$.  If $A \le 0, \operatorname{Rank} (B)=n$ and
$\int_{\{|z-z_0|\le\varepsilon\}} \rho_0(z)^{-1}\,\mathrm{d}z<\infty$ holds for
some $z_0\in \mathbb{R}^d$ and some $\varepsilon>0$, then the associated Markov
transition probability $P_t(x,\dsd{y})$ satisfies
\[
\|P_t (x, \cdot)- P_t (y, \cdot)\|_{\mathrm{var}} \le
\frac{C(1+|x-y|)}{\sss t},\qquad x,y\in \mathbb{R}^d, t>0,
\]
for some constant $C>0$, which is sharp for large $t$ and implies
that  the process has successful couplings. The Harnack inequality,
ultracontractivity and the strong Feller property are also investigated
for the (conditional) transition semigroup.
\end{abstract}


\end{frontmatter}

\section{Introduction}\label{sec1}

L\'{e}vy processes are fundamental models of Markov processes, from
which more general diffusion-jump-type Markov processes can be
constructed by solving stochastic differential equations or martingale
problems. It is well known that a L\'{e}vy process can be decomposed
into two independent parts, that is, the Brownian (or Gaussian) part
and the jump part. Comparing with the analysis on the Brownian motion,
that on the pure jump part is far from complete. For instance, except
for stable-like processes that can be treated as subordinations of
diffusion processes \cite{SV} (see also \cite{BJ,JS} for heat kernel
upper bounds for $\aal$-stable processes with drifts),  little is known
concerning regularities  of the transition probabilities of O--U-type
jump processes. Most existing regularity results for O--U (or
generalized Mehler) semigroups were derived by   using the Gaussian
part as the leading term (cf. \cite{DZ,ROW,RW03} and references
within). In contrast, besides known results on the transition density
for L\'{e}vy processes (see \cite{KS,FV,Tu} and references therein),
the strong Feller property was recently proved by Priola and Zabczyk
\cite{PZ} for O--U jump processes  by considering a H\"{o}rmander
condition and L\'{e}vy measures. The main purpose of this paper is to
investigate more regular properties on O--U semigroups in the same
spirit,  so that our results work well for the pure jump case as
emphasized in the Abstract.

Recall that a L\'{e}vy measure   $\nu$ on $\R^d$ is such that
$\nu(\{0\})=0$ and (see \cite{A})
\[
\int_{\R^d} (|z|^2\land 1)\nu(\dsd{z})<\infty.
\]
Let $b\in \R^d$ and $Q$
be a non-negatively definite $d\times d$ matrix. The underlying
L\'{e}vy process $L_t$ is the Markov process on $\R^d$ generated by
\[
\scr L f := \langle b,\nn f\rangle + \operatorname{Tr}(Q\nn^2 f)+\int_{\R^d} \bigl\{f(z+\cdot)-f
-\langle \nn f, z\rangle 1_{\{|z|\le 1\}}\bigr\}\nu(\dsd{z}),
\]
which is well
defined for $f\in C_b^2(\R^d)$.

Now, let $A$ be a real  $n\times n$  matrix and $B$ be a real $n\times
d$ matrix. We shall investigate the solution to the following linear
stochastic differential equation
\begin{equation}\label{1.0}
\dsd{X}_t^x= (A X_t^x)\,\dsd{t} + B\,\dsd{L}_t,\qquad X_0^x=x\in\R^n.
\end{equation}
We shall investigate the following properties of
the solution:
\begin{enumerate}[(A)]
\item[(A)] The coupling property.
\item[(B)] The Harnack inequality and ultracontractivity.
\item[(C)] The strong Feller property.
\end{enumerate}

The coupling method is a powerful tool in the study of Markov
processes, and the coupling property that we are going to study is
closely related to long-time behaviors, Liouville-type properties and
the 0--1 law of tail-/shift-invariant events. Recall that a couple
$(X_t,Y_t)$ is called a coupling of the Markov process associated with
a given transition probability if both $X_t$ and $Y_t$ are  Markov
processes  associated with the same transition probability (possibly
with different initial distributions). In this case, $X_t$ and $Y_t$
are called the marginal processes of the coupling.  A~coupling
$(X_t,Y_t)$ is called successful if the coupling time
\[
\mathbf T:=\inf\{t\ge 0\dvtx X_t=Y_t\}<\infty,\qquad\mbox{a.s.}
\]
A Markov  process is said to have a coupling property (or to have
successful couplings) if, for any initial distributions $\mu_1$ and
$\mu_2,$  there exists  a successful coupling with marginal processes
starting from  $\mu_1$ and $\mu_2$, respectively.   A slightly weaker
notion is the shift-coupling property: for any two initial
distributions there exists a coupling $(X_t,Y_t)$ with marginal
processes  starting from them respectively, such that ``$X_{\mathbf
T_1}= Y_{\mathbf T_2}$'' holds for some finite stopping times $\mathbf
T_1, \mathbf T_2$ (see \cite{AT}). In general, the coupling property is
stronger than the shift-coupling property, but they are equivalent  if
the Markov semigroup  satisfies a weak parabolic Harnack  inequality
(see \cite{CW}).

Consider a strong Markov process with transition semigroup $P_t$. For
any (not necessarily successful) coupling with initial distributions
$\mu_1$ and $\mu_2$, one has  (see \cite{L,Chen})
\begin{equation}\label{1.1"}
\|\mu_1P_t -\mu_2P_t\|_{\mathrm{var}}\le 2\ppbl(\mathbf{T}>t),\qquad t\ge 0,
\end{equation}
where $\|\cdot\|_{\mathrm{var}}$ is the
total variational norm. This follows by setting $X_t=Y_t$ for $t\ge
\mathbf T$ due to the strong Markov property. Moreover, for any
coupling $(X_t,Y_t)$ with initial distributions $\dd_x$ and $\dd_y$, a
bounded harmonic function $f$ (i.e. $P_t f=f$ for $t\ge 0$) satisfies
\[
|f(x)-f(y)|\le \inf_{t>0} \E|f(X_t)-f(Y_t)|\le 2\|f\|_\infty \ppbl(\mathbf T=\infty).
\]
Consequently,
if a strong Markov process  has a coupling property, then its bounded
harmonic functions have to be constant, that is, the Liouville property
holds for   bounded harmonic functions. In general, the coupling
property of a
 strong Markov process on
$\R^n$ with semigroup $P_t$ is equivalent to each of the following
statements (see  \cite{CG}, Section \ref{sec4}, and \cite{L}, Chapters 3 and 5):
\begin{enumerate}[(iii)]
\item[(i)] For any $\mu_1,\mu_2\in \scr P(\R^n),\lim_{t\to\infty} \|\mu_1 P_t-\mu_2 P_t\|_{\mathrm{var}}=0.$
\item[(ii)] All bounded time--space harmonic functions are constant, that is, a bounded measurable function $u$ on
$[0,\infty)\times \R^n$ has to be constant if $u(t,\cdot)= P_s u(t+s,\cdot)$ holds for all $s,t\ge 0.$
\item[(iii)] The tail $\si$-algebra of the process is trivial, that is,
$\ppbl(X\in A)=0\mbox{ or }1$ holds for any initial distribution and any
$A\in \bigcap_{t>0} \si((\R^n)^{[0,\infty)}\ni  w\mapsto  w_s\dvtx s\ge t).$
\end{enumerate}
Correspondingly,  each of the following statements are
equivalent to the shift-coupling property (see \cite{T}, Section \ref{sec4}, or
\cite{AT}):
\begin{enumerate}[(vi)]
\item[(iv)] For any $\mu_1,\mu_2\in \scr P(\R^n),\lim_{t\to\infty}
\|\ff 1 t \int_0^t(\mu_1-\mu_2) P_s\,\dsd{s}\|_{\mathrm{var}}=0$.
\item[(v)] All bounded harmonic functions are constant.
\item[(vi)] The shift-invariant  $\si$-algebra of the process is trivial, that is,
$\ppbl(X \in B)=0\mbox{ or }1$ for any initial distribution and any
shift-invariant measurable set $B\subset (\R^n)^{[0,\infty)}$.
\end{enumerate}
In Section \ref{sec3}, we shall present explicit conditions on $A, B$ and the
L\'{e}vy measure $\nu$ such that the coupling property holds (see
Theorem \ref{T3.1}).

Next, we aim to establish the following Harnack inequality for $P_t$
initiated in \cite{W97} for diffusion semigroups:
\[
(P_t f(x))^\aal\le (P_t f^\aal(y)) H_\aal(t, x,y),\qquad t>0, x,y\in \R^n,
\aal>1,
\]
for positive measurable functions $f$, where $H_\aal$ is a positive function on $(0,\infty)\times (\R^n)^2.$

When $\nu(\R^d)<\infty$, with a positive probability the process does
not jump before a fixed time $t>0$, so that
 this inequality could not hold for the pure jump case. This is the main reason why all existing results in this
 direction only work for the case with a non-degenerate Gaussian part
 (cf. \cite{ROW,RW03}). To work out the Harnack inequality also for the pure jump case, we shall be restricted on
 the event that the process jumps before time $t$. More precisely, let $\tau_1$ be the first jump time of  the L\'{e}vy
 process induced by an absolutely continuous part of $\nu$. If $\Rank (B)=n$, then the Harnack inequality and
 ultracontractivity are investigated in Section \ref{sec4} for the following modified
sub-Markov operator $P_t^1$ (see Theorem \ref{T4.1}):
\[
P_t^1 f(x):= \E\bigl\{ f(X_t^x) 1_{\{t\ge \tau_1\}}\bigr\}.
\]

Finally, we look at the strong Feller property of $P_t$. By the same
reasoning that leads to the invalidity of the Harnack inequality, when
$\nu$ is finite the pure jump semigroup cannot be strong Feller.
Therefore, in \cite{PZ} the authors only considered the case that $\nu$
is infinite. More precisely, if $\nu$ has an infinite absolutely
continuous part and if there exists $m\ge 1$ such that the rank
condition
\[
\Rank(B, AB,\ldots, A^{m-1}B)=n
\]
holds, then \cite{PZ}, Theorem 1.1 and Proposition 2.1, imply the
strong Feller property of $P_t$.  We shall extend this result by
allowing the absolutely continuous part of $\nu$ to be finite. In this
case the number $m$ in the rank condition will refer to the strong
Feller property of the semigroup conditioned by the event that the
$m$th jump happens before time $t>0$ (see Theorem \ref{T5.1}).

It might be interesting to indicate that for jump processes the strong
Feller property is incomparable with the coupling property. Indeed, the
latter is a long-time property but the former is somehow a short-time
property. For the strong Feller property, we need the process to be
able to visit any area before any fixed time, for which  the jump
measure has to be infinite as  mentioned above. However, the situation
in the diffusion case is very different: Whenever the diffusion process
is able to visit any area for a long time, it will be able to do so
before any fixed time.

The remainder of the paper is organized as follows. To study the
coupling property and the Harnack inequality, we shall first
investigate in Section \ref{sec2} the quasi-invariance of random shifts for
compound Poisson processes, which in particular leads to a conditional Girsanov theorem. Then we will study the properties included in (A), (B)
and (C) in Sections \ref{sec3}, \ref{sec4} and \ref{sec5}, respectively.

\section{Quasi-invariance and the Girsanov theorem}\label{sec2}

Throughout of this section, we assume $\llam:=\nu(\R^d)\in (0,\infty)$
and let $L:=\{L_t\}_{t\ge 0}$ be the compound
 Poisson process with L\'{e}vy measure $\nu$ and $L_0=0$. Let $\LL$ be the distribution of $L$,
which is a probability measure on the path space
\[
W=\Biggl\{\sum_{i=1}^\infty x_i1_{[t_i,\infty)}\dvtx i\in \N, x_i\in \R^d\setminus\{0\}, 0\le t_i\uparrow \infty\mbox{ as }i\uparrow \infty\Biggr\}
\]
equipped with the $\si$-algebra induced by $\{ w\mapsto  w_t\dvtx t\ge 0\}$.
Let $\DD  w_t=  w_t- w_{t-}$ for $t>0$.

For any $T>0$, let $\LL_T$ be the distribution of $L_{[0,T]}
:=\{L_t\}_{t\in [0,T]}$, and let $\tau$ and $\xi$ be random variables
with distributions $\ff 1 T 1_{[0,T]}(t)\,\dsd{t}$
on $[0,T]$ and $\ff 1 \llam \nu$ on $\R^d$, respectively, such that $\xi,\tau$ and $L$ are independent. It is shown
in \cite{WY} that the distribution of $L_{[0,T]}+ \xi 1_{[\tau, T]}$ is $\ff 1 {\llam T} n_T (w)\LL_T(\dsd{w})$, where
\begin{equation}\label{2.0}
n_T( w):= \#\{t\in [0,T]\dvtx w_t\ne  w_{t-}\},\qquad w\in W.
\end{equation}
We shall extend this result to more general random
variables $\xi$ and $\tau$. To this end, write

\begin{equation}\label{1.1}
L_t= \sum_{i=1}^{N_t} \xi_i,\qquad t\ge 0,
\end{equation}
where $N_t$ is the Poisson process on $\Z_+$ with rate $\llam$ and
$\{\xi_i\}_{i\ge 1}$ are i.i.d. random variables on $\R^d$, which are
independent of $N:=\{N_t\}_{t\ge 0}$ and have common distribution $\ff
1\llam \nu$.

\begin{tthm}\label{T2.1}
Let $(\xi,\tau)$ be a random variable on
$\R^d\times [0,\infty)$. Then the distribution of
$L+\xi1_{[\tau,\infty)}$ is absolutely continuous with respect to $\LL$
if and only if the joint distribution of $(L, \xi,\tau)$ has the form
\[
\vv \LL(\dsd{w}) \dd_0(\dsd{z}) \Theta( w,\dsd{t})+ g( w, z,
t)\LL(\dsd{w})\nu(\dsd{z})\,\dsd{t},
\]
where $\vv\in [0,1]$ is a constant, $g$
is a non-negative measurable function on $W\times \R^d\times
[0,\infty)$ and $\Theta( w,\dsd{t})$ is a transition probability from $W$
to $[0,\infty).$ In this case,  the distribution of $L+\xi
1_{[\tau,\infty)}$ is formulated as
\[
\biggl\{\vv + \sum_{\DD  w_t\ne 0} g\bigl( w-\DD w_t1_{[t,\infty)},\DD w_t, t\bigr)\biggr\}\LL(\dsd{w}).
\]
\end{tthm}

According to Theorem \ref{T2.1}, the random shift $L\mapsto
L+\xi1_{[\tau,\infty)}$ is quasi-invariant if and only if the
conditional distribution of $(L,\xi,\tau)$ given $\{\xi\ne 0\}$ is
absolutely continuous w.r.t. the product measure
$\LL\times\nu\times\dsd{t}$. Since, when $\xi=0$, the random shift does
not help the coupling, below we will only choose non-zero $\xi$.

To prove this result, we shall make use of the Mecke formula for the
Poisson measure.  Let $E$ be a Polish space with Borel $\si$-field
$\F$, and let $\si$ be a locally finite measure on $E$. Then $\pi_\si$,
the Poisson measure with intensity $\si$,  is a probability measure on
the configuration space
\[
\GG:= \Biggl\{\sum_{i=1}^n\dd_{x_i}\dvtx n\in \Z_+\cup\{\infty\}, x_i\in E\Biggr\}
\]
fixed by the Laplace transform
\[
\int_{\GG} \e^{\ggam(f)}\pi_\si(\dsd\ggam) = \e^{\si(\e^f-1)},\qquad f\in C_0(E).
\]
Note that the corresponding $\si$-field on $\GG$ is induced by
$\{\ggam\mapsto \ggam(f)\dvtx f\in C_0(E)\}.$ The Mecke formula \cite{M} (see
also \cite{R}) says that for any non-negative measurable function $F$
on $\GG\times E$,
\begin{equation}\label{M}
\int_{\GG\times E}F(\ggam+\dd_z,z)\pi_\si(\dsd\ggam)\si(\dsd{z})
=\int_\GG \pi_\si(\dsd\ggam)\int_E F(\ggam,z)\ggam(\dsd{z}).
\end{equation}

By considering the Poisson measure with intensity $\nu\times \dsd{t}$ on
$\R^d\times [0,\infty)$, we will be able to prove Theorem \ref{T2.1} by
using the following result.

 \begin{tthm}\label{TP} Let
$A\subset E$ be measurable, $X$ be a random variable on $\GG$ with
distribution $\pi_\si$ and $\eta$ be a random variable on $E$. Then the
measure $\ppbl(X+\dd_\eta\in\cdot, \eta\in A)$ is absolutely continuous
with respect to $\pi_\si$ if and only if the  measure
 $\ppbl((X,\eta)\in\cdot, \eta\in A)$ is absolutely continuous with respect to
$\pi_\si\times\si$.
\end{tthm}

\begin{pf}
Let $D_{X,\eta}$ be the distribution of $(X,\eta)$.

(a) The sufficiency.  Assume that $\ppbl((X,\eta)\in \cdot, \eta\in A)=
g(\ggam,z)\pi_\si(\dsd\ggam)\si(\dsd{z})$ for some non-negative measurable
function $g$ on $\GG\times E$. For any bounded measurable function $f$
on $\GG$, by the Mecke formula (\ref{M}) for
\[
F(\ggam,z):= f(\ggam) g(\ggam-\dd_z,z)1_{\{\ggam\ge \dd_z\}}1_A(z),
\]
we have
\begin{eqnarray*}
\E \{1_A(\eta) f(X+\dd_\eta)\}
&=&
\int_{\GG\times A}f(\ggam+\dd_z)g(\ggam,z)\pi_\si(\dsd\ggam)\si(\dsd{z})
\\
&=&
\int_\GG f(\ggam) \biggl\{\int_A g(\ggam-\dd_z,z)\ggam(\dsd{z})\biggr\}\pi_\si(\dsd\ggam).
\end{eqnarray*}
So,
$\ppbl(X+\dd_\eta\in\cdot, \eta\in A)$ is absolutely continuous with
respect to  $\pi_\si$ with density function  $\ggam\mapsto \int_A
g(\ggam-\dd_z,z)\ggam(\dsd{z}).$

(b) The necessity. Assume that $\ppbl(X+\dd_\eta\in\cdot, \eta\in A)$ is
absolutely continuous with respect to $\pi_\si$. For any measurable set
$N\subset \GG\times E$ with $(\pi_\si\times\si)(N)=0,$ we intend to
prove
\begin{equation}\label{AB}
\ppbl\bigl((X,\eta)\in N, \eta\in A\bigr)=D_{X,\eta}(N)=0.
\end{equation}
Let
\begin{eqnarray*}
&&A_N= \{\ggam+\dd_z\dvtx (\ggam,z)\in N, z\in A\}\subset \GG,\qquad F(\ggam,z)=
1_{N\cap (\GG\times A)}(\ggam-\dd_z,z)\\
&&\quad\mbox{for }(\ggam,z)\in \GG\times E.
\end{eqnarray*}
If $\ggam\in A_N$, then there exists $z_0\in A$ such that
$(\ggam-\dd_{z_0},z_0)\in N.$ This means that $\ggam\ge \dd_{z_0}$ and
\[
\int_EF(\ggam,z)\ggam(\dsd{z})\ge h(\ggam,z_0)=1.
\]
Therefore,
\[
\int_EF(\ggam,z)\ggam(\dsd{z})\ge 1_{A_N}(\ggam),\qquad\ggam\in \GG.
\]
Combining this with (\ref{M}) and noting that
$(\pi_\si\times\si)(N)=0$, we obtain
\begin{eqnarray*}
\pi_\si(A_N)&\le&
\int_\GG\pi_\si(\dsd\ggam)\int_E F(\ggam,z)\ggam(\dsd{z}) =\int_{\GG\times E}
F(\ggam+\dd_z, z) \pi_\si(\dsd\ggam)\si(\dsd{z})\\
&\le&\int_{\GG\times E} 1_N(\ggam,z)
\pi_\si(\dsd\ggam)\si(\dsd{z})=(\pi_\si\times\si)(N)=0.
\end{eqnarray*}
Since $\ppbl(X+\dd_\eta\in\cdot, \eta\in A)$ is absolutely continuous with
respect to $\pi_\si$, this implies that
\[
\ppbl\bigl((X,\eta)\in N, \eta\in A\bigr)\le \ppbl(X+\dd_\eta\in A_N,\eta\in A)=0.
\]
Thus, (\ref{AB}) holds.
\end{pf}

\begin{pf*}{Proof of Theorem \ref{T2.1}}
(1) The sufficiency. Let
$\pi_\si$ be the Poisson measure with intensity $\si:= \nu(\dsd{z})\times
\dsd{t}$. Since $\si(\{0\})=0$   and the Lebesgue measure $\dsd{t}$ is
infinite on $[0,\infty)$ without atom, $\pi_\si$ is supported
on
\[
\GG_0:= \Biggl\{\sum_{i=1}^\infty \dd_{(x_i,t_i)}\dvtx i\in \N, x_i\in \R^d\setminus\{0\},
0\le t_i\uparrow \infty\mbox{ as }i\uparrow \infty\Biggr\}.
\]
Let
\[
\psi\dvtx W\to \GG_0;\qquad\sum_{i=1}^\infty x_i 1_{[t_i,\infty)} \mapsto \sum_{i=1}^\infty \dd_{(x_i,t_i)}.
\]
We have (see  \cite{B}, page 12)
\begin{equation}\label{2.1}
\pi_\si= \LL\circ\psi^{-1},\qquad\LL=\pi_\si\circ\psi.
\end{equation}
By (\ref{M}), for any non-negative measurable
function $h$ on
$\GG_0\times \R^d\times [0,\infty)$,  we have
\[
\int_{\GG_0}\pi_\si(\dsd\ggam)\int_{\R^d\times [0,\infty)} h(\ggam,x,t) \ggam(\dsd{x}, \dsd{t}) =
\int_{\GG_0\times \R^d\times [0,\infty)} h\bigl(\ggam+ \dd_{(x,t)}, x, t\bigr)
\pi_\si(\dsd \ggam) \nu(\mathrm{d} x)\,\dsd{t}.
\]
Combining this with (\ref{2.1}) we
conclude that
\[
\int_W \sum_{\DD  w_t\ne 0} H( w, \DD w_t, t)\LL(\dsd{w}) = \int_{W\times \R^d\times [0,\infty)}
H\bigl( w+ x1_{[t,\infty)}, x, t\bigr) \LL(\dsd{w})\nu(\dsd{x})\,\dsd{t}
\]
holds for any
non-negative measurable function $H$ on $W\times \R^d\times
[0,\infty)$. Therefore, for
 any non-negative measurable function $F$ on $W$,
we have
\begin{eqnarray*}
&&
\E F\bigl(L+\xi 1_{[\tau,\infty)}\bigr)
\\
&&\quad=
\E\bigl\{F(L)1_{\{\xi=0\}}\bigr\}+\int_{W\times\R^d\times [0,\infty)} F\bigl( w + x1_{[t,\infty)}\bigr) g( w, x,t)\LL(\dsd{w}) \nu(\dsd{x})\,\dsd{t}
\\
&&\quad=
\int_WF( w) \biggl\{\vv + \sum_{\DD w_t\ne 0} g\bigl( w-\DD w_t 1_{[t,\infty)},\DD w_t, t\bigr)\biggr\}\LL(\dsd{w}).
\end{eqnarray*}
This completes the
proof of the sufficiency.

(2) The necessity. Let the distribution of $L+\xi1_{[\tau,\infty)}$ be
absolutely continuous with respect to $\LL$. Let $\vv=\ppbl(\xi=0)$ and
let $\Theta( w,\dsd{t})$ be the regular conditional distribution of
$\tau$ given $L$ and $\xi=0$. Then for any non-negative measurable
function $f$ on $W\times \R^d\times [0,\infty)$,
\[
\E f(L,\xi,\tau)=\vv \int_{\GG\times [0,\infty)} f( w, 0, t)
\LL(\dsd{w})\Theta( w, \dsd{t})+\E \bigl\{f(L,\xi,\tau)1_{\{\xi\ne 0\}}\bigr\}.
\]
So,
to prove that the distribution of $(L,\xi,\tau)$ has the required form,
it suffices to show that for any $\LL\times\nu\times\dsd{t}$-null set
$N$, we have
\begin{equation}\label{2}
\ppbl\bigl((L,\xi,\tau)\in N, \xi\ne 0\bigr)=0.
\end{equation}
To this
end, we shall make use of Theorem \ref{TP}. Let $E= \R^d\times
[0,\infty)$ and $ X=\psi(L)$. We have
\[
\psi\bigl(L+\xi1_{[\tau,\infty)}\bigr)=X+\dd_{(\xi,\tau)}\qquad\mbox{for } \xi\ne 0.
\]
Let
\[
\ttil N=\{(\psi( w),z,t)\dvtx ( w,z,t)\in N, z\ne 0\}.
\]
By (\ref{2.1}) we have
\begin{equation}\label{LLL}
(\pi_\si\times \nu\times\dsd{t})(\ttil N)\le(\LL\times\nu\times\dsd{t})(N)=0.
\end{equation}
Now, since the
distribution of $L+\xi1_{[\tau,\infty)}$ is absolutely continuous with
respect to $\LL$, due to (\ref{2.1}) so is
$\ppbl(X+\dd_{(\xi,\tau)}\in\cdot,$ $\xi\ne 0)$ with respect to $\pi_\si$.
Hence, according to Theorem \ref{TP}, $\ppbl((X,\xi,\tau)\in\cdot,$ $\xi\ne
0)$ is absolutely continuous with respect to $\pi_\si\times
\nu\times\dsd{t}$. Combining this with (\ref{2.1}) and (\ref{LLL}), we
arrive at
\[
\ppbl\bigl((L,\xi,\tau)\in N, \xi\ne 0\bigr)=\ppbl\bigl((X,\xi,\tau)\in \ttil N, \xi\ne 0\bigr)=0.
\]
Therefore, (\ref{2}) holds.
\end{pf*}

In the situation of Theorem \ref{T2.1}, let
\begin{equation}\label{U}
U(w)= \vv+ \sum_{\DD w_t\ne 0} g\bigl( w-\DD w_t 1_{[t,\infty)}, \DD w_t,t\bigr),\qquad w\in W.
\end{equation}
As a direct consequence of Theorem \ref{T2.1},
the following result says that the distribution of $L+\xi
1_{[\tau,\infty)}$ under probability
\[
\ff{1_{\{U>0\}}}{\LL(U>0)U}\bigl(L+\xi1_{[\tau,\infty)}\bigr) \ppbl
\]
coincides with that of $L$ under probability
$\ff{1_{\{U(L)>0\}}}{\LL(U>0)}\ppbl$. This can be regarded as  a
conditional Girsanov theorem.

\begin{cor} \label{C2.2}
In the situation of Theorem $\ref{T2.1}$  let
$U$ be in $(\ref{U})$. Then for any non-negative measurable function
$F$ on $W$,
\[
\E\bigl\{\bigl(F1_{\{U>0\}}\bigr)(L)\bigr\}
= \E\biggl\{\ff{F1_{\{U>0\}}}{U}\bigl(L+\xi1_{[\tau,\infty)}\bigr) \biggr\}.
\]
\end{cor}

\section{The coupling property}\label{sec3}

Recall that for the Brownian motion the equality in  (\ref{1.1"}) is
reached by the coupling by reflection covered by Lindvall and Rogers in
\cite{LR}. More precisely, let $P_t^B(x,\dsd{y}) $ be the transition
probability of the Brownian motion on $\R^d$ and let $\mathbf T_{x,y}$
be the coupling time of the coupling by reflection for initial
distributions $\dd_x$ and $\dd_y$. One has (see \cite{CL}, Section~5)
\begin{eqnarray}\label{B}
\ff 1 2 \|P_t^B(x,\cdot) - P_t^B(y,\cdot)\|_{\mathrm{var}}
&=&
\ppbl(\mathbf T_{x,y}>t)\nonumber
\\
&=&
\ff{\sss 2}{\sss \uppi} \int_0^{|x-y|/(2\sss t)}\e^{-u^2/2}\,\dsd{u}\le \ff{\sss 2|x-y|}{\sss t},\qquad t>0.\quad
\end{eqnarray}
Our first result aims to provide an analogous estimate for L\'{e}vy
jump processes, which in particular implies the coupling property of
the process according to the equivalent statement (i). Intuitively, to
ensure the coupling property for a  L\'{e}vy jump process, the L\'{e}vy
measure should have a non-discrete support to make the process active
enough. In this paper, we shall assume that $\Rank (B) =n$ and $\nu$
has a non-trivial absolutely continuous part.

\begin{tthm}\label{T3.1}
Let $P_t(x,\dsd{y})$  be the transition probability
for the solution to $(\ref{1.0})$. Let  $\Rank (B)=n$ and $\langle Ax,
x\rangle \le 0$ hold for $x\in \R^n$. If $\nu\ge \rr_0(z)\,\dsd{z}$ such
that
\[
\int_{\{|z-z_0|\le\vv\}}\rr_0(z)^{-1}\,\dsd{z}<\infty
\]
holds for some $z_0\in\R^d$ and some $\vv>0$, then
\begin{equation}\label{LL}
\|P_t (x,\cdot)- P_t (y,\cdot)\|_{\mathrm{var}} \le
\ff{C(1+|x-y|)}{\sss t},\qquad x,y\in \R^n, t>0,
\end{equation}
holds for some
constant $C>0$, and hence, the coupling property and assertions \textup{(i)--(vi)} hold.
\end{tthm}

\begin{remm}
(1) According to
\cite{PZ04}, Theorem 3.5(ii), if $A$ has an eigenvalue with a positive
real part, then, under an assumption on large jumps, the coupling
property fails. In this sense the assumption $A\le 0$ is somehow
reasonable for the coupling property. On the other hand, by
\cite{PZ04}, Theorem 3.8, in the diffusion case, all bounded harmonic
functions could be constant (i.e., the shift-coupling property holds)
provided all eigenvalues of $A$ have non-positive real parts. It
would be interesting to extend this result to the jump case.

(2)\vspace*{1pt} The condition $\int_{\{|z-z_0|\le\vv\}}\rr_0(z)^{-1}\,\dsd{z}<\infty$
follows from $\inf_{|z-z_0|\le\vv} \rr_0(z)>0$, which corresponds to
the uniformly elliptic condition in the diffusion setting. Similarly to
(\ref{B}) in the Brownian
 motion case, (\ref{LL})  is   sharp for large $t>0$ in the pure jump case. To see this, let $n=d=B=1, A=Q=0$
 and let $\nu$ be a probability measure such that
\[
\int_\R z \nu(\dsd{z})=0,\qquad\int_{\R} z^2\nu(\dsd{z})=1,\qquad\int_{\R} |z|^3\nu(\dsd{z})<\infty.
\]
Then the corresponding L\'{e}vy process reduces to the compound Poisson
process up to a constant drift $b_0$:
\[
X_t= \sum_{i=1}^{N_t}\xi_i +b_0t,
\]
where $N_t$ is the Poisson process on $\Z_+$ with rate $1$ and
$\{\xi_i\}_{i\ge 1}$ are i.i.d. and independent of $N_t$ with common
distribution $\nu$. By the Berry--Esseen inequality (see \cite{S}),
\[
\sup_{r\in \R} \bigl|\ppbl\bigl( X_t< r\sss t + b_0t\bigr) -\Phi(r)\bigr|\le \ff {c_0}{\sss t},\qquad
t>0,
\]
holds
for some constant $c_0>0$, where $\Phi$ is the standard Gaussian
distribution function. Therefore,
\begin{eqnarray*}
\|P_t(x,\cdot)-P_t(0,\cdot)\|_{\mathrm{var}}
&\ge&
2 \sup_{r\in \R}\bigl|\ppbl\bigl(X_t <r\sss t+b_0t\bigr)- \ppbl\bigl( X_t <r\sss t + b_0t -x\bigr)\bigr|
\\
&\ge&
2\sup_{r\in\R} \bigl|\Phi(r)-\Phi\bigl(r-x/\sss t\bigr)\bigr| -\ff{4 c_0} {\sss t}\ge\ff{c_1|x|-4 c_0}{\sss t},\qquad t\ge
x^2,
\end{eqnarray*}
holds for some constant $c_1>0.$
\end{remm}

It is well known that the solution to (\ref{1.0}) can be formulated as

\begin{equation}\label{3.0}
X_t^x = \e^{At} x + \int_0^t \e^{A(t-s)} B\,\dsd{L}_s,\qquad x\in \R^n, t\ge 0.
\end{equation}
To make use of Theorem \ref{T2.1}, we shall split $L_t$ into two
independent parts:
\[
L_t=L_t^1+L_t^0,
\]
where $L^0:=\{L_t^0\}_{t\ge 0}$ is the compound Poisson process with L\'{e}vy measure
$\nu_0(\dsd{z}):= \rr_0(z)\,\dsd{z}$, and $L^1:=\{L_t^1\}_{t\ge 0}$ is the
L\'{e}vy  process with L\'{e}vy measure $\nu-\nu_0$ generated by $\scr
L-\scr L_0$ for
\[
\scr L_0 f:= \int_{\R^d} \bigl(f(\cdot +z)-f(z)\bigr)\nu_0(\dsd{z}).
\]
So, (\ref{3.0}) reduces to
\begin{equation}\label{3.0'}
X_t^x = \e^{At} x + \int_0^t \e^{A(t-s)} B\,\dsd{L}_s^1 +
\int_0^t \e^{A(t-s)} B\,\dsd{L}_s^0,\qquad x\in \R^n, t\ge 0.
\end{equation}
Moreover, let
\begin{equation}\label{1.1'}
L_t^0 = \sum_{i=1}^{N_t}\xi_i,
\end{equation}
where $N:=\{N_t\}_{t\ge 0}$ is the Poisson process
on $\Z_+$ with rate $\llam_0:=\nu_0(\R^d),$ and $\{\xi_i\}_{i\ge 1}$ are
i.i.d. real random variables with common distribution $\nu_0/\llam_0$
such that $N, \{\xi_i\}_{i\ge 1}$ and $L^1$ are independent.

To prove Theorem \ref{T3.1},  we  introduce the following fundamental
lemma.

\begin{lem}\label{L3.1}
Let  $\llam_0\in (0,\infty)$, and let
$\{\eta_i\}_{i\ge 1}$ be a sequence of square-integrable real random
variables that are conditional independent given $N$ such that
 $\E(\eta_i|N)=1$ and $\E(\eta_i^2|N)\le \si$ hold for some constant $\si\in (0,\infty)$  and all $i\ge 1$. Then
\[
E\Biggl(1-\ff 1 {\llam_0 T} \sum_{i=1}^{N_T} \eta_i\Biggr)^2 \le \ff \si {\llam_0 T}.
\]
\end{lem}

\begin{pf}
Since $\E(\eta_i\eta_j|N)=1$ for $i\ne j$ and
$\E(\eta_i^2|N)\le \si$ for  $i\ge 1$, we have
\begin{eqnarray*}
\E \Biggl(1-\ff 1 {\llam_0 T}\sum_{i=1}^{N_T}\eta_i\Biggr)^2
&=&
\ff 1 {(\llam_0T)^2} \E\Biggl(\sum_{i=1}^{N_T}\eta_i\Biggr)^2 -\ff 2{\llam_0T}\E\sum_{i=1}^{N_T} \eta_i +1
\\
&=&
\ff 1 {(\llam_0T)^2} \E\Biggl\{ \sum_{i,j=1}^{N_T} \E(\eta_i\eta_j|N)\Biggr\} -\ff 2 {\llam_0T}\E\Biggl\{\sum_{i=1}^{N_T} \E(\eta_i|N)\Biggr\}+1
\\
&\le&
\ff 1 {(\llam_0T)^2} \sum_{n=1}^\infty \ff{(n^2-n+\si n)(\llam_0 T)^n\e^{-\llam_0T}}{n!} -1 = \ff{\si}{\llam_0T}.
\end{eqnarray*}
\upqed\end{pf}

\begin{pf*}{Proof of Theorem \ref{T3.1}}
We simply denote $\B_r=\{z\dvtx |z-z_0|\le r\}$ for $r>0.$ Using $\rr_0\land 1$ to replace
$\rr_0$, we may and do assume that $\rr_0\le 1.$ In this case
$\nu_0(\dsd{z}):=\rr_0(z)\,\dsd{z}$ is finite. For $T>0$, let $\tau$ be a random variable on $[0,\infty)$
with distribution $\ff 1 T 1_{[0,T]}(t)\,\dsd{t}$ and $\xi$ on $\R^n$ with distribution
\[
\ff{1_{\B_{\vv/2}}(z)\nu_0(\dsd{z})}{\nu_0(\B_{\vv/2})},
\]
such that $L^0, L^1,\xi,\tau$ are independent. Let $\LL(\dsd{w})$ be
the distribution of $L^0$. It is easy to see that the distribution of
$(L^0, \xi,\tau)$ is
\[
\ff {1_{\B_{\vv/2}}(z)1_{[0,T]}(t)} {T\nu_0(\B_{\vv/2})}\LL(\dsd{w})\nu_0(\dsd{z})\,\dsd{t}.
\]
By Theorem \ref{T2.1} and (\ref{1.1'}), for any $z\in\R^d$  we  have
\begin{eqnarray*}
&&
\E f\biggl(\e^{AT}z  +\int_0^T\e^{A(T-t)}B\,\dsd \bigl(L^0+\xi 1_{[\tau,\infty)}\bigr)_t\biggl)
\\
&&\quad=
\E\biggl\{\ff{f(\e^{AT}z  +\int_0^T\e^{A(T-t)}B\,\dsd{L}^0_t)}{T\nu_0(\B_{\vv/2})}\sum_{t\le T} 1_{\B_{\vv/2}\setminus\{0\}}(\DD L_t^0)\biggr\}
\\
&&\quad=
\E\biggl\{\ff{f(\e^{AT}z  +\int_0^T \e^{A(T-t)}B\,\dsd{L}^0_t)}{T\nu_0(\B_{\vv/2})} \sum_{i=1}^{N_T}1_{\B_{\vv/2}}(\xi_i)\biggr\}.
\end{eqnarray*}
Letting $\pi_{x,T}$
be the distribution of $x+\int_0^T\e^{-At}B\,\dsd{L}_t^1$ and combining
this with (\ref{3.0'}) and the independence of $L^0$ and $L^1$, we
obtain
\begin{eqnarray}\label{3.1}
&&\E f\bigl(X_T^x+\e^{A(T-\tau)}B\xi\bigr)\nonumber
\\
&&\quad=
\E f\biggl(\e^{AT}\biggl\{x+ \int_0^T\e^{-At}B\,\dsd{L}_t^1\biggr\}+\int_0^T\e^{A(T-t)}B\,\dsd \bigl(L^0+\xi1_{[\tau,\infty)}\bigr)_t\biggr)\nonumber
\\
&&\quad=
\int_{\R^d}\biggl\{\E f\biggl(\e^{AT}z  +\int_0^T \e^{A(T-t)}B\,\dsd\bigl(L^0+ \xi 1_{[\tau,\infty)}\bigr)_t\biggr)\biggr\}\pi_{x,T}(\dsd{z})\nonumber
\\
&&\quad=
\int_{\R^d}\E \Biggl\{\ff{f(\e^{AT}z  +\int_0^T\e^{A(T-t)}B\,\dsd{L}^0_t)}{T\nu_0(\B_{\vv/2})} \sum_{i=1}^{N_T}1_{\B_{\vv/2}}(\xi_i)\Biggr\}\pi_{x,T}(\dsd{z})\nonumber
\\
&&\quad=
\E \Biggl\{\ff{f(X_T^x )}{T\nu_0(\B_{\vv/2})} \sum_{i=1}^{N_T}1_{\B_{\vv/2}}(\xi_i)\Biggr\}.
\end{eqnarray}
Next, since $\Rank (B)=n$, we have $d\ge n$ and up to a permutation of
coordinates in $\R^d$, we may and do assume that $B=(B_1, B_2)$ for
some invertible $n\times n$   matrix $B_1$ and some $n\times (d-n)$
matrix $B_2.$ If, in particular, $n=d$, then $B_1=B$. Moreover, for
simplicity we write
\[
\R^n= \R^n\times \{\bar 0\}\subset \R^d,
\]
where $\bar 0$ is the original in $\R^{d-n}$. In other words,
for any $x\in \R^n$, we set $x=(x,\bar 0)\in \R^d.$ Since $\langle
Ax,x\rangle \le 0$ for $x\in \R^n$, if
\[
\|B_1^{-1}\|\cdot |x-y|\le \ff \vv 2,
\]
then
\[
|B_1^{-1}\e^{\tau A}(x-y)|\le \|B_1^{-1}\|\cdot |\e^{\tau A}(x-y)|\le \|B_1^{-1}\|\cdot |x-y|\le \ff \vv 2.
\]
So    the distribution of $(L^0, \xi +B_1^{-1}\e^{A\tau} (x-y), \tau )$
is

\begin{eqnarray*}
&&
\ff {1_{[0,T]}(t) 1_{\B_{\vv/2} +B_1^{-1}\e^{At}(x-y)}(z)}{T\nu_0(\B_{\vv/2})}\LL(\dsd{w}) \nu_0\bigl(\dsd{z}- B_1^{-1}\e^{At}(x-y)\bigr)\,\dsd{t}
\\
&&\quad=
\ff { 1_{[0,T]}(t)1_{\B_{\vv/2}+B_1^{-1}\e^{At}(x-y)}(z)\rr_0(z-B_1^{-1}\e^{At}(x-y))}{T\nu_0(\B_{\vv/2})\rr_0(z)}\LL(\dsd{w})\nu_0(\dsd{z})\,\dsd{t}.
\end{eqnarray*}
Similarly to (\ref{3.1}), due to  Theorem \ref{T2.1}, (\ref{1.1'}) and
the independence of $L^0$ and $L^1$, we have
\begin{eqnarray*}
&&
\E f\bigl(X_T^x+\e^{A(T-\tau)}B\xi\bigr)
\\
&&\quad=
\E f\biggl(\e^{AT}y + \int_0^T \e^{A(T-t)}B\,\dsd \bigl(L^1+L^0+ \{\xi +B_1^{-1}\e^{A\tau}(x-y)\}1_{[\tau,\infty)}\bigr)_t\biggr)
\\
&&\quad=
\E\biggl\{\ff{f(X_T^y)}{T\nu_0(\B_{\vv/2})} \sum_{t\le T}1_{(\B_{\vv/2}+B_1^{-1}\e^{At}(x-y))\setminus\{0\}}
(\DD L_t^0)\ff{\rr_0(\DD L_t^0-B_1^{-1}\e^{At}(x-y))}{\rr_0(\DD L_t^0)}\biggr\}
\\
&&\quad=
\E\Biggl\{\ff{f(X_T^y)}{T \nu_0(\B_{\vv/2})} \sum_{i=1}^{N_T}1_{\B_{\vv/2}+B_1^{-1}\e^{A\tau_i}(x-y)}
(\xi_i)\ff{\rr_0(\xi_i-B_1^{-1}\e^{A\tau_i}(x-y))}{\rr_0(\xi_i)}\Biggr\},
\end{eqnarray*}
where $\tau_i$ is the $i$th jump time of $N_t$ for $i\ge 1$. Combining
this with (\ref{3.1}), we arrive at
\begin{eqnarray}\label{3.2}
&&
|P_Tf(x)- P_Tf(y)|\quad\nonumber
\\
&&\quad\le
\E\Biggl| f(X_T^y) \Biggl(1- \ff 1 {T \nu_0(\B_{\vv/2}) } \sum_{i=1}^{N_T}1_{\B_{\vv/2}+B_1^{-1}\e^{A\tau_i}(x-y)}(\xi_i)\ff{\rr_0(\xi_i-B_1^{-1}\e^{A\tau_i} (x-y))}{\rr_0(\xi_i)}\Biggr)
\quad\nonumber\\
&&\qquad\hphantom{\E|}{}+ f(X_T^x)\Biggl(\ff 1 {T\nu_0(\B_{\vv/2})} \sum_{i=1}^{N_T}
1_{\B_{\vv/2}}(\xi_i) -1\Biggr) \Biggr|,\qquad |x-y|\le \ff \vv
{2\|B_1^{-1}\|}.\quad
\end{eqnarray}
To apply Lemma \ref{L3.1}, let
\[
\eta_i= \ff{\llam_0 1_{\B_{\vv/2}}(\xi_i)}{\nu_0(\B_{\vv/2})},\qquad \ttil\eta_i
= \ff{\llam_0
\rr_0(\xi_i-B_1^{-1}\e^{A\tau_i}(x-y))}{\nu_0(\B_{\vv/2})\rr_0(\xi_i)}
1_{\B_{\vv/2}+B_1^{-1}\e^{A\tau_i}(x-y)}(\xi_i),i\ge 1.
\]
Then
$\{\eta_i\}_{i\ge 1}$ are i.i.d. and independent of $N$ with
\begin{eqnarray*}
\E \eta_i
&=&
\ff 1 {\llam_0} \int_{\B_{\vv/2}}\ff{\llam_0}{\nu_0(\B_{\vv/2})} \nu_0(\dsd{z})=1,
\\
\E \eta_i^2
&=&
\ff 1 {\llam_0}\int_{\B_{\vv/2}}\ff{\llam_0^2}{\nu_0(\B_{\vv/2})^2} \nu_0(\dsd{z})=\ff{\llam_0}{\nu_0(\B_{\vv/2})}<\infty,
\end{eqnarray*}
while $\{\ttil\eta_i\}_{i\ge 1}$ are conditional independent given $N$ such that
\begin{eqnarray*}
\E(\ttil\eta_i|N)
&=&
\ff 1 {\llam_0}\int_{\B_{\vv/2}+B_1^{-1}\e^{A\tau_i}(x-y)}\ff{\llam_0\rr_0(z-B_1^{-1}\e^{A\tau_i}(x-y))}{\nu_0(\B_{\vv/2})\rr_0(z)}\nu_0(\dsd{z})
\\
&=&
\ff 1 {\llam_0} \int_{\B_{\vv/2}+B_1^{-1}\e^{A\tau_i}(x-y)}\ff{\llam_0\rr_0(z-B_1^{-1}\e^{A\tau_i}(x-y))}{\nu_0(\B_{\vv/2})}\,\dsd{z}=1,
\end{eqnarray*}
and since $\rr_0\le 1$ and $|B_1^{-1}\e^{A\tau_i}(x-y)|\le \ff \vv
2$,
\begin{eqnarray*}
\E(\ttil\eta_i^2|N)
&=&
\ff 1 {\llam_0}\int_{\B_{\vv/2}+B_1^{-1}\e^{A\tau_i}(x-y)}\ff{\llam_0^2\rr_0(z-B_1^{-1}\e^{A\tau_i}(x-y))^2}{\nu_0(\B_{\vv/2})^2\rr_0(z)^2}\nu_0(\dsd{z})
\\
&\le&
\llam_0\int_{\B_{\vv/2}+B_1^{-1}\e^{A\tau_i}(x-y)}\ff{\dsd{z}}{\nu_0(\B_{\vv/2})^2\rr_0(z)}\\
&\le&
\ff{\llam_0}{\nu_0(\B_{\vv/2})^2} \int_{\B_\vv}\ff{\dsd{z}}{\rr_0(z)}<\infty.
\end{eqnarray*}
Therefore, by (\ref{3.2}) and Lemma \ref{L3.1},
\[
\|P_T(x,\cdot)-P_T(y,\cdot)\|_{\mathrm{var}}\le \ff c {\sss T},\qquad T>0, |x-y|\le \ff\vv
{2\|B_1^{-1}\|},
\]
holds for some constant $c>0.$ This implies (\ref{LL}) for some constant $C>0$
since for $|x-y|>\ff\vv {2\|B_1^{-1}\|} $ and $m_{x,y}:= \inf\{i\in \N\dvtx i\ge 2\|B_1^{-1}\|\cdot |x-y|/\vv\},$ we have
\begin{eqnarray*}
&&
\|P_t(x,\cdot)-P_t(y,\cdot)\|_{\mathrm{var}}
\\
&&\quad\le
\sum_{i=1}^{m_{x,y}}\biggl\|P_t\biggl(x+ \ff{i(y-x)}{m_{x,y}},\cdot\biggr)-P_t\biggl(x+\ff{(i-1)(y-x)}{m_{x,y}},\cdot\biggr)\biggr\|_{\mathrm{var}}.
\end{eqnarray*}

Finally, it is easy to see that (\ref{LL}) implies the statement (i)
and hence, the coupling property of the process.
\end{pf*}

To conclude this section, we present a result on the equivalence of the
coupling property and the shift-coupling property by using a criterion
in \cite{CW}.

\begin{prp} Let $\nu(\R^d)<\infty$ and $A=0$. If either $b=\int_{\{|z|\le
1\}} z\nu(\dsd{z})$ or $\Rank (B)=n$ and $Q$ is non-degenerate, then
the coupling property is equivalent to the shift-coupling property.
\end{prp}

\begin{pf}
Let $\llam:=\nu(\R^d)<\infty$ and $A=0$.  Let
$L_t=L_t^1+L_t^0$ as before for $L_t^0$ being  the compound Poisson
process   specified in (\ref{1.1'}) for $\nu$ in place of $\nu_0$.  If
$Q=0$ and
 $b= \int_{\{|z|\le 1\}}z\nu(\dsd{z})$, then $X_t^x= x+BL_t^0$. So, for any non-negative measurable function
 $f$ on $\R^n$, and  any $t,s>0$, we have
\begin{equation}\label{L1}
P_{t +s}f(x)= \E f(BL^0_{t+s} + x) \ge \E
\bigl[f(BL^0_t+x)1_{\{N_{t+s}-N_t=0\}}\bigr] = \e^{-\llam s}P_tf(x).
\end{equation}
Therefore, by  \cite{CW}, Theorem 5, the coupling property is
equivalent to the shift-coupling property.

Next, let $A=0$, $Q$ be non-degenerate and $\Rank (B)=n$. Let $P_t^J$
and $P_t^D$ be the semigroups of $B  L_t^0$ and $BL_t^1$, respectively.
Then it is easy to see that the generator of $P_t^D$ is an elliptic
operator with constant coefficients and hence, satisfies the
Bakry--Emery curvature-dimension condition. Therefore, according to
\cite{BQ}, there exists a constant $k\ge n$ such that
\begin{equation}\label{X}
P_t^Df\le \biggl(\ff{t+s} s\biggr)^{k/2} P_{t+s}^Df,\qquad t,s>0,
\end{equation}
holds for non-negative measurable function $f$.  Since $A=0$ implies
that the diffusion part and the jump part are independent, we have
$P_t= P_t^DP_t^J,$ where $P_t^J$ is the semigroup associated to
$BL_t^0$, which satisfies (\ref{L1}). Therefore,
\[
P_{t} f\le \biggl(\ff{t+s} s\biggr)^{k/2}\e^{\llam s} P_{t+s}f,\qquad
f\ge 0, s,t>0.
\]
This implies the equivalence of the coupling property and the
shift-coupling property according to \cite{CW}, Theorem 5.
\end{pf}

The condition $b=\int_{\{|z|\le 1\}} z\nu(\dsd{z})$ is used to ensure
the desired inequality (\ref{L1}). If this condition does not hold,
there exists $b_0\ne 0$ such that $X_t^x= x +BL_t^0 +b_0t$, so that
instead of (\ref{L1}) one has
\[
P_{t +s}f(x)\ge\E\bigl[f(BL^0_t+x+b_0s)1_{\{N_{t+s}-N_t=0\}}\bigr]= \e^{-\llam s}P_tf(x+b_0s),
\]
which is not enough  to apply \cite{CW}, Theorem 5.

\section{Harnack inequality and ultracontractivity}\label{sec4}

Let $\nu\ge \nu_0:= \rr_0(z)\,\dsd{z}>0$ for some $\rr_0>0$ with
$\llam_0:=\nu_0(\R^d)\in (0,\infty).$ As in Section \ref{sec3}, let
$L_t=L_t^1+L_t^0$ such that
 $L^0$ and $L^1$ are independent, where  $L^0$ is the
compound Poisson process with L\'{e}vy measure $\nu_0.$   Let  $\tau_1$
be the first jump time of $L_t^0$. We shall establish the Harnack
inequality for
\begin{equation}\label{N}
P_t^1 f(x):= \E\bigl\{f(X_t^x)1_{\{\tau_1\le t\}}\bigr\},\qquad t\ge 0, x\in \R^n,f\in \scr B_b(\R^n).
\end{equation}

\begin{tthm}\label{T4.1}
Let $\nu\ge \nu_0:= \rr_0(z)\,\dsd{z}$ with
$\llam_0:=\nu_0(\R^d)\in (0,\infty)$, and let $P_t^1$ be defined above.
Let $ \Rank (B)=n$. There exists a constant $c=c(B)>0$ such that
if
\[
V_p(r):= \ff 1 {\llam_0}\sup_{|z'|\le r} \int_{\R^d}\ff{
\rr_0(z-z')^{p/(p-1)}}{\rr_0(z)^{1/(p-1)} }\,\dsd{z} <\infty,\qquad r\ge
0,
\]
holds for some $p>1$, then for any positive measurable function $f$ on
$\R^n$,
\[
( P_t^1 f(x))^p \le ( P_t^1 f^p(y))\bigl\{(1-\e^{-\llam_0 t})V_p\bigl(c\e^{\|A\|t}|x-y|\bigr)\bigr\}^{p-1},\qquad x,y\in \R^d,
t>0,
\]
holds. Consequently,
\[
\|P_t^1\|_{p\to \infty} \le(1-\e^{-\llam_0 t})\e^{\|A\|t/p}
\biggl\{\int_{\R^d}\ff {\dsd{x}} {V_p(c\e^{\|A\| t}|x|)^{p-1}}
\biggr\}^{-1/p}<\infty,
\]
where $\|\cdot\|_{p\to q}$ is the operator norm
from $L^p(\R^n;\dsd{x})$ to $L^q(\R^n;\dsd{x})$ for any $p,q\ge
1.$
\end{tthm}

\begin{pf}
Let $L^0, L^1,\xi,\tau$ be independent such that the
distributions of $\xi$ and $\tau$ are $\nu_0/\llam_0$
and $\ff 1 T 1_{[0,T]}(t)\,\dsd{t}$, respectively.  As in the proof of Theorem \ref{T3.1}, let $B=(B_1, B_2)$
such that $B_1$ is invertible.  Since the distribution of $(L^0, \xi, \tau)$ is
\[
\ff{1_{[0,T]}(t)}{\llam_0 T} \LL(\dsd{w}) \nu_0(\dsd{z})\,\dsd{t},
\]
Corollary \ref{C2.2} holds for
\[
U( w)= \ff 1 {\llam_0 T} n_T( w),
\]
where $n_T$ is defined by (\ref{2.0}). Since $\tau\le T$ and $\xi\ne 0,$ which
are independent of $L^0$ and $L^1$, we have
\[
U\bigl(L^0+\xi 1_{[\tau,\infty)}\bigr) =\ff 1 {\llam_0 T} n_T\bigl(L^0+ \xi1_{[\tau,\infty)}\bigr)>0.
\]
Therefore, by Corollary \ref{C2.2}
and noting that $\tau_1\le T$ a.s. for the process
$L^0+\xi1_{[\tau,\infty)}$,
\begin{eqnarray}\label{3.6}
P_T^1 f(x)
&=&
\E\bigl[f(X_T^x) 1_{\{\tau_1\le T\}}\bigr]\nonumber
\\
&=&
\E\biggl\{\ff{\llam_0 T f(\e^{AT}x+\int_0^T\e^{A(T-t)}B\,\dsd (L^1+L^0 +\xi1_{[\tau,\infty)})_t)}{n_T(L^0+\xi 1_{[\tau,\infty)})}\biggr\}\nonumber
\\
&=&
\E\biggl\{\ff{\llam_0 T f(\e^{AT}y+\int_0^T\e^{A(T-t)}B
\,\dsd(L^1+L^0 +\{\xi+B_1^{-1}\e^{A\tau}(x-y)\}1_{[\tau,\infty)})_t)}{n_T(L^0+\{\xi+B_1^{-1}\e^{A\tau}(x-y)\}1_{[\tau,\infty)})}\biggr\}
\nonumber\\
&=&
\E\Biggl\{\ff{f(X_T^y)1_{\{\tau_1\le T\}}}{N_T} \sum_{i=1}^{N_T}\ff{\rr_0(\xi_i-B_1^{-1} \e^{A\tau_i}(x-y)}{\rr_0(\xi_i)}\Biggr\},
\end{eqnarray}
 where  $c=\|B_1^{-1}\|.$ By the H\"{o}lder inequality, we obtain
\begin{eqnarray*}
(P_T^1 f(x))^p
&\le&
P_T^1 f^p(y)\Biggl\{\E\Biggl( \ff{1_{\{N_T\ge 1\}}}{N_T} \sum_{i=1}^{N_T}\ff{\rr_0(\xi_i-B_1^{-1}\e^{A\tau_i}(x-y))}{\rr_0(\xi_i)}\Biggr)^{p/(p-1)}\Biggr\}^{p-1}
\\
&\le&
P_T^1 f^p(y)\Biggl\{\sum_{n=1}^\infty\ff{(\llam_0T)^n \e^{-\llam_0 T}}{n (n!)}\sum_{i=1}^n\sup_{|z'|\le c\e^{\|A\|T}|x-y|}\E\biggl(\ff{\rr_0(\xi_i-z')}{\rr_0(\xi_i)}\biggr)^{p/(p-1)}\Biggr\}^{p-1}
\\
&=&
P_T^1 f^p(y)\bigl\{(1-\e^{-\llam_0 T})V_p\bigl(c\e^{\|A\|T}|x-y|\bigr)\bigr\}^{p-1}.
\end{eqnarray*}
This implies the desired Harnack inequality.

Next, since there exists a probability $\mu_T$ on $\R^n$ such that
\[
P_T f^p(x)=: \E f^p(X_T^x)=\int_{\R^n} f^p(\e^{AT} x+y) \mu_T(\dsd{y}),
\]
if $\int_{\R^n} f^p(x)\,\dsd{x}\le 1,$ then
\[
\int_{\R^n} P_T^1 f^p(x)\,\dsd{x} \le \int_{\R^n} P_T f^p(x)\,\dsd{x} =\int_{\R^n}\mu_T(\dsd{y})
\int_{\R^n} f^p(\e^{TA}x+y)\,\dsd{x}\le \e^{\|A\|T}.
\]
Therefore, by the Harnack inequality, for any non-negative $f$ with
$\int_{\R^d} f^p(z)\,\dsd{z}\le 1$,
\begin{eqnarray*}
&&
(P_T^1 f(x))^p \int_{\R^d} \ff {\dsd{y}} {(V_p(c\e^{\|A\|T}|x-y|))^{p-1}}
\\
&&\quad\le
(1-\e^{-\llam_0 T})^{p-1} \int_{\R^d} P_T^1 f^p(y) \,\dsd{y} \le (
1-\e^{-\llam_0 T})^p\e^{\|A\|T}.
\end{eqnarray*}
This implies
the desired upper bound of $\|P_T^1\|_{p\to\infty}.$
\end{pf}

It is easy to see that $V_p<\infty$ holds for many concrete choices of
$\rr_0$, including $\rr_0(z):= c_1 \e^{-c_2 |z|^r}$ for some constants
$c_1,c_2,r>0$  and $\rr_0(z):= c (1+|z|)^{-r}$ for some $r>d$ and
$c>0.$

Finally, when $\nu$ has a large enough absolutely continuous part, we
may derive the ultracontractivity by comparing with the $\alpha$-stable
process.

\begin{tthm}
Assume that $n=d$ and $B=I$. Let $\aal\in (0,2).$ If
\[
\nu(\dsd{z})\ge \ff{c}{|z|^{\aal+d}}1_{\{|z|<r\}}\,\dsd{z}
\]
holds for some
constants $c,r>0$, then
\[
\|P_t\|_{1\to\infty}\le \ff{c'}{(1\land t)^{d/\aal}},\qquad t>0,
\]
holds for some constant $c'>0.$
\end{tthm}

\begin{pf}
(a) We first observe that if $r=\infty$, that is,
\begin{equation}\label{AL}
\nu(\dsd{z})\ge \ff{c}{|z|^{\aal+d}}\,\dsd{z},
\end{equation}
then
\[
\|P_t\|_{1\to\infty}\le \ff{c'}{t^{\dsd/\aal}},\qquad t\in (0,1],
\]
holds. When $A=0$ and $\nu(\dsd{z})\ge \ff{c}{|z|^{\aal+d}}\,\dsd{z}$ this is
well known according to the heat kernel upper bound of the
$\alpha$-stable process. In general, let $\eta$ be the symbol of the
L\'{e}vy process $L$ with characteristics $(b,Q,\nu)$. Let $\mu_t$ be
the probability measure on $\R^d$ with Fourier transform
\[
\hat \mu_t(z)= \exp\biggl[-\int_0^t \eta(\e^{sA^*}z)\,\dsd{s}\biggr],\qquad z\in\R^d.
\]
We have
\[
P_tf(x)= \int_{\R^d} f(\e^{tA}x +y) \mu_t(\dsd{y}).
\]
Let  $c_1>0$ be such that
\[
c\int_0^t|\e^{s A^*} z|^\aal \,\dsd{s}\ge c_1t |z|^\aal,\qquad t\in [0,1].
\]
According to (\ref{AL}) there are two probability measures $\mu_t^1$
and $\mu_t^2$ on $\R^d$ such that $\mu_t=\mu_t^1*\mu_t^2$ and the
Fourier transform of $\mu_t^1$ is
\[
\hat \mu_t^1(z)=\exp[-c_1t|z|^\aal].
\]
Combining this with
the known heat kernel bound  of the $\alpha$-stable process, we can
find a constant $c'>0$ such that for any $f\ge 0$,
\begin{eqnarray*}
P_t f(x)
&=&
\int_{\R^d} \mu_t^1(\dsd{z})\int_{\R^d} f(\e^{tA} x+ y+ z)\mu_t^2(\dsd{y})
\\
&\le&
\ff {c'}{t^{d/\aal}} \int_{\R^d} f(z)\,\dsd{z},\qquad x\in \R^d,t\in
(0,1].
\end{eqnarray*}
This implies the desired estimate.

(b) Let $r\in (0,\infty)$. To apply (a), let $L^0$ be the compound
Poisson process independent of $L$ with L\'{e}vy measure
\[
\nu_0(\dsd{z}):= \ff{c}{(|z|\lor r)^{d+\aal}}\,\dsd{z}.
\]
Then $\bar L:=
L+L^0$ is a L\'{e}vy process with L\'{e}vy measure
\[
\bar\nu(\dsd{z})= \nu(\dsd{z})+\nu_0(\dsd{z})\ge \ff{c}{|z|^{\aal+d}}\,\dsd{z}.
\]
Let $\bar P_t$ be the semigroup associated with the equation
\[
\dsd \bar X_t= A\bar X_t\,\dsd{t} +\dsd \bar L_t.
\]
By (a)
\begin{equation}\label{ABC}
\|\bar P_t\|_{1\to\infty} \le \ff{c'}{t^{d/\aal}},\qquad t\in (0,1],
\end{equation}
holds for some constant $c'>0.$ Let $\tau_1$
be the first jump time of $L^0$. We have
\begin{eqnarray*}
\bar P_t f(x)
&:=&
\E f\biggl(\e^{At}x +\int_0^t \e^{A(t-s)}\,\dsd{L}_s +\int_0^t \e^{A(t-s)}\,\dsd{L}_s^0\biggr)
\\
&\ge&
\E \biggl\{1_{\{\tau_1>t\}}f\biggl(\e^{At}x + \int_0^t\e^{A(t-s)}\,\dsd{L}_s\biggr)\biggr\}
\\
&=&
\e^{-\llam_0 t} P_t f(x),\qquad f\ge 0,
\end{eqnarray*}
where
$\llam_0:=\nu_0(\R^d)<\infty$. Combining this with (\ref{ABC}) we
complete the proof.
\end{pf}

\section{Strong Feller property}\label{sec5}

As in Sections \ref{sec3} and \ref{sec4}, let $\nu\ge \nu_0:=\rr_0(z)\,\dsd{z}$ for some
non-negative measurable function $\rr_0$ on $\R^d$ such that
$\llam_0:=\nu_0(\R^d)>0$. Let $L_t= L_t^1+L_t^0$   for independent
$L^1$ and $L^0$ such that $L_t^0$ is the compound Poisson process with
L\'{e}vy measure $\nu_0$. For any $i\ge 1$, let $\tau_i$ be the $i$th
jump time of $L_t^0$. If $\llam_0=\infty$, we set $\tau_i=0$ for all
$i\ge 1$ by convention. We shall prove the strong Feller property for
the operator  $P_{t }^{m}$ defined by
\begin{equation}\label{5.1}
P_{t }^m f(x)= \E\bigl\{f(X_t^x)1_{\{  \tau_{m}\le t\land(\tau_1+
t_m)\}}\bigr\},
\end{equation}
where $m\ge 1$ and
\[
t_m:=\sup\{t\ge 0:\Rank(\e^{s_1 A}B,\ldots, \e^{s_m A}B)=n,\forall 0\le s_1<\cdots<s_m\le t\}.
\]
According to the following lemma, we have $t_m>0$ provided the rank
condition
\begin{equation}\label{H}
\Rank (B, AB, \ldots, A^{m-1}B)=n
\end{equation}
holds. This extends \cite{PZ}, Lemma 2.2, by allowing $m\ne n$.

\begin{lem} \label{L5.2}
If $(\ref{H})$ holds for some $m\ge 1$, then
$t_m>0$. Consequently, for $0\le s_1<\cdots <s_m\le t_m$ and
\[
\psi_{s_1,\ldots, s_m}(z_1,\ldots,z_m):= \sum_{i=1}^m \e^{s_i
A}Bz_i,\qquad z_1,\ldots, z_m\in \R^d,
\]
$\ggam\circ\psi_{s_1,\ldots, s_m}^{-1}$
is an absolutely continuous probability measure on $\R^n$ provided so
is $\ggam$ on $\R^{md}$.
\end{lem}

\begin{pf}
By \cite{PZ}, Lemma 2.3, it suffices to prove the first
assertion. For $0\le s_1<\cdots<s_m$, let
\begin{eqnarray*}
F_{i,0}^{(0)}
&=&
\e^{s_iA},\qquad 1\le i\le m,
\\
F_{i, k}^{(k)}
&=&
\ff{F_{i, k-1}^{(k-1)}-F_{k,k-1}^{(k-1)}}{s_i-s_k},\qquad 1\le k\le m-1,k+1\le i\le m.
\end{eqnarray*}
Since
\[
\ff{\dsd^i}{\dsd{s}^i}\e^{sA} \biggl|_{s=0} =  A^i,\qquad i\ge 0,
\]
for any $1\le i \le m,$  $F_{i, i-1}^{(i-1)}$
approximates $A^{(i-1)}$ as $s_m\downarrow 0.$ Therefore, there exist real matrices $U_1,\ldots, U_m$
depending on $(s_1,\ldots, s_m)$ such that
\[
\lim_{s_m\to 0} \|U_i\|=0,\qquad 1\le i\le m,
\]
and
\[
F_{i, i-1}^{(i-1)} = A^{i-1} +U_i,\qquad 1\le i\le m.
\]
Since $\{F_{i,i-1}^{(i-1)}\dvtx 1\le i\le m\}$
are linear combinations of
$\{\e^{s_i A}\dvtx 1\le i\le m\},$ we have
\begin{equation}\label{R}
\Rank(\e^{s_1A}B, \e^{s_2 A}B,\ldots, \e^{s_mA}B)\ge
\Rank(B+ U_1B, AB+U_2B, \ldots, A^{m-1}B +U_mB).
\end{equation}
Since
$(B, AB, \ldots, A^{m-1}B)$ has full rank $n$, and since $U_iB\to 0$ as
$s_m\to 0$, there exists $t>0$ such that if $0\le s_1<\cdots<s_m\le t$,
then
\[
\Rank(B+U_1B, AB+U_2B,\ldots, A^{m-1}B+U_mB)=n.
\]
Combining this with (\ref{R}) we complete the proof.
\end{pf}

\begin{tthm}\label{T5.1}
If $t_m>0$, then  $P_{t }^{m}$ is strong Feller
for $t>0.$ Consequently, if $(\ref{H})$ holds for some $m\ge 1$, then
 $t_{m\land n}>0$ such that $P_t^{n\land m}$ is strong
Feller for $t>0.$
\end{tthm}

\begin{pf}
According to Lemma \ref{L5.2}
and the fact that (\ref{H}) with $m\ge n$ is equivalent to the
condition with $m=n$ (cf. \cite{Z}), it suffices to prove the first
assertion. We shall complete the proof in four easy steps.

(a) We first observe that $P_{t }^m$ is strong Feller if
\[
P_{t }^m(0,\dsd{x}):= \ppbl(X_t^0\in \dsd{x},t\ge \tau_m, t_m\ge \tau_m-\tau_1)
\]
is absolutely continuous.
Indeed, let $P_{t }^m(0,\dsd{x})= g(x)\,\dsd{x}.$ Then
\[
P_{t }^m f(x)= \int_{\R^n}f(\e^{At}x +y) g(y)\,\dsd{y}.
\]
Therefore, $P_{t }^m$ is strong Feller according to \cite{H}, Lemma 11.

(b) Next, we claim that it suffices to prove the result for $\llam_0<\infty.$ If $\llam_0=\infty,$ then for any
$l\ge 1$ let $\nu_l= (\rr_0\land l)(z)\,\dsd{z}$ and $\llam_l= \nu_l(\R^d).$ Let $\tau_i(l)$ be the $i$th jump time
for the corresponding compound Poisson process with L\'{e}vy measure $\nu_l$. If the assertion holds for finite
$\llam_0$, then we may use $\nu_l$ to replace $\nu_0$ so that
\[
\ppbl\bigl(X_t^0\in\dsd{x}, t\ge \tau_m(l), t_m\ge \tau_m(l)-\tau_1(l)\bigr)
\]
is absolutely continuous. Therefore,
for any measurable set $D\subset \R^n$ with volume $|D|=0$,
\[
P_t(0,D)\le \ppbl\bigl(X_t^0\in D, t\ge \tau_m(l), t_m\ge \tau_m(l)-\tau_1(l)\bigr)+ \ppbl\bigl(\tau_m(l)\ge t\land t_m\bigr)= \e^{-\llam_l(t\land t_m)/m}.
\]
Since $\llam_l\uparrow \llam_0=\infty$ as $l\uparrow \infty$, we see that $P_t(0,\cdot)$
is absolutely continuous.

(c) We aim to show that it suffices to prove for the case that
$L_t=L_t^0$, that is, $\nu=\nu_0$ and the
 L\'{e}vy process is the compound Poisson process with L\'{e}vy measure $\nu_0.$
Indeed, since
\[
X_t^0= \int_0^t \e^{(t-s)A}B\,\dsd{L}_s^1+ \int_0^t\e^{(t-s)A}B\,\dsd{L}_s^0,
\]
where  $L^1$ and $L^0$ are independent,
$P_{t }^m(0,\dsd{x})$ is absolutely continuous provided so is
\[
\ppbl\biggl( \int_0^t\e^{(t-s)A}B\,\dsd{L}_s^0\in \dsd{x}, t\ge \tau_m, t_m\ge \tau_m-\tau_1\biggr).
\]

(d) Now, assume that  $\nu=\nu_0$ with $\llam_0\in (0,\infty)$ and
$L_t=L_t^0$. Let $\pi(\dsd{s}_1,\ldots,\dsd{s}_m)$ be
the distribution of $(\tau_1,\ldots,\tau_m),$ and let
\[
K=\{(s_1,\ldots, s_m)\dvtx s_m-s_1\le t_m, 0< s_1<\cdots s_m\le t\}.
\]
Since by (\ref{3.0'}) and (\ref{1.1'}) with $L_t^1=0$
\[
X_t^0= \int_0^t\e^{(t-s)A}B\,\dsd{L}_s^0= \e^{(t-\tau_m)A}\sum_{i=1}^m\e^{(\tau_m-\tau_i)A}B\xi_i +\int_{\tau_m}^t\e^{(t-s)A}B\,\dsd{L}_s^0
\]
provided $\tau_m\le t$, for any non-negative measurable function $f$ on $\R^n$,
we have
\begin{equation}\label{AA}
P_{t }^m f(0) = \int_K \E
f\Biggl(\e^{(t-s_m)A}\sum_{i=1}^m \e^{(s_m-s_i)A}B\xi_i +\int_{s_m}^t
\e^{(t-s)A}B\,\dsd{L}_s^0\biggr)\pi(\dsd{s}_1,\ldots,\dsd{s}_m),
\end{equation}
where
$\{\xi_i\}$ are i.i.d. random variables with distribution $\nu_0/\llam_0$
independent of $(L_s^0)_{s\ge s_m}$. Since $\e^{(t-s_m)A}$ is
invertible and $s_m-s_i< t_m$, by the definition of $t_m$ the mapping
\[
(z_1,\ldots, z_m)\mapsto \e^{(t-s_m)A}\sum_{i=1}^m \e^{(s_m-s_i)A}Bz_i
\]
is onto, so that the distribution
of the random variable
\[
\e^{(t-s_m)A}\sum_{i=1}^m \e^{(s_m-s_i)A}B\xi_i
\]
is absolutely continuous (see \cite{PZ}, Lemma 2.3). By (\ref{AA}) and the independence
of this random variable and
\[
\int_{s_m}^t
\e^{(t-s)A}B\,\dsd{L}_s^0,
\]
we conclude that $P_{t }^m(0,\dsd{x})$ is
absolutely continuous.
\end{pf}

\begin{remm}
In concrete examples we may have
$t_m=\infty$ so that $P_t^m$ reduces to
\[
P_t^m f(x):= \E\bigl\{f(X_t^x)1_{\{\tau_{m}\le t\}}\bigr\},
\]
which refers to the conditional distribution of $X_t^x$
in the event that $L_t^0$ jumps at least $m$ times before $t$. For
instance, as in \cite{PZ}, formula (1.3), let $n=2, d=1$ and
\[
A= \pmatrix{ 0  &1\cr1 &0},\qquad B=\pmatrix{0\cr 1}.
\]
We have $A^2 =I$ and $AB = {{0}\choose{1}}.$ So,
\[
\e^{sA}B = \Biggl(\sum_{n=0}^\infty \ff{s^{2n}}{(2n)!} \Biggr) \pmatrix{ 0\cr
1}
+ \Biggl(\sum_{n=0}^\infty \ff{s^{2n+1}}{(2n+1)!}\Biggr)\pmatrix{
1\cr 0}=\cosh(s)\pmatrix{ 0\cr
1}+ \sinh(s) \pmatrix{ 1\cr
0}
\]
 holds for all $s\ge 0.$ Since $\sinh(s_2-s_1)>0$ for  $s_2-s_1>0$, and since $\e^{s_1A}$ is invertible, we have
\[
\Rank(\e^{s_1 A}B, \e^{s_2 A}B) = \Rank \bigl(B, \e^{(s_2-s_1)A}B\bigr)= 2=n.
\]
 Therefore, $t_2=\infty$.
 \end{remm}

 \section*{Acknowledgements}
 The author would like to thank Enrico Priola and a referee for  helpful comments
  on the first version of the paper and for introducing him to the very interesting paper \cite{PZ}. This paper supported in
 part by WIMCS and SRFDP.

\printhistory

\end{document}